# Integer Sequences associated with Integer Monic Polynomial

**Ashok Kumar Gupta**
Department of Electronics and Communication,
Allahabad University, Allahabad - 211 002, India
(Email address: akgjkiapt@hotmail.com)

**Ashok Kumar Mittal**
Department of Physics,
Allahabad University, Allahabad – 211 002, India
(Email address: mittal_a@vsnl.com)

**Abstract:** To every integer monic polynomial of degree m can be associated m integer sequences having interesting properties to the roots of the polynomial. These sequences can be used to find the real roots of any integer monic polynomial by using recursion relation involving integers only. This method is faster than the conventional methods using floating point arithmetic.
.

## 1. Introduction

A polynomial is said to be integer monic [1] if all its coefficients are integers and the coefficient of highest power is unity. To a given polynomial can be associated an *m*-term recursion relation for generating integer sequences. A set of *m* such sequences, which together exhibit interesting properties related to the roots of the polynomial, can be obtained if the *m* initial terms of each of these *m* sequences is chosen in a special way using the companion matrix [2] of the polynomial. The companion matrix of a polynomial is a matrix such that the characteristic equation of the companion matrix is the given polynomial.

## 2. Construction of the integer sequences

Let $p(x) = x^m + a_1 x^{m-1} + \ldots + a_m$, where $a_i$ are integers, be the given monic polynomial. The companion matrix of the polynomial is given by

$$R = \begin{bmatrix} -a_1 & -a_2 & -a_3 & \ldots & -a_{m-1} & -a_m \\ 1 & 0 & 0 & \ldots & 0 & 0 \\ 0 & 1 & 0 & \ldots & 0 & 0 \\ . & . & . & \ldots & . & . \\ . & . & . & \ldots & 0 & . \\ 0 & 0 & 0 & \ldots & 1 & 0 \end{bmatrix} \qquad (1)$$

Let $r_i$, $i = 1, 2, \ldots, m$ be the roots of $p(x)$. Then each $r_i$ is an eigenvalue of **R** with corresponding eigenvector $[r_i^{m-1}, r_i^{m-2}, \ldots, r_i, 1]^T$. The matrix **R** satisfies the equation

$$\mathbf{R}^m + a_1 \mathbf{R}^{m-1} + \ldots + a_m = \mathbf{0} \tag{2}$$

Let

$$\mathbf{S}_0 = (S^{(1)}_0, S^{(2)}_0, \ldots, S^{(m)}_0)^T \tag{3}$$

where $S^{(i)}_0$ are arbitrary integers. Define

$$\mathbf{S}_j = \mathbf{R}^j \mathbf{S}_0 \ , \ j = 1, 2, \ldots, m-1 \tag{4}$$

The first $m$ terms of the sequence $S^{(i)}$ are taken to be $\{S^{(i)}_0, S^{(i)}_1, \ldots, S^{(i)}_{m-1}\}$. Beyond this the $m$-term recurrence relation defines the terms of the sequence

$$S^{(i)}_j = -a_1 S^{(i)}_{j-1} - a_2 S^{(i)}_{j-2} - \ldots - a_m S^{(i)}_{j-m} \ , \ j > m - 1, \ i = 1, 2, \ldots, m \tag{5}$$

These sequences satisfy the interesting property that each of the ratios $S^{(i)}_j$ to $S^{(i+1)}_j$, $i = 1, 2, \ldots, m-1$ tends to the root of $p$, having the largest absolute value, if it is unique.

## 3. Example

Let

$$p(x) = x^2 + 2x - 1 \tag{6}$$

be the given integer monic polynomial. Here $m = 2$, $a_1 = 2$ and $a_2 = -1$. The companion matrix to the polynomial (6) is given by

$$\mathbf{R} = \begin{bmatrix} -2 & 1 \\ 1 & 0 \end{bmatrix} \tag{7}$$

We take

$$\mathbf{S}_0 = [1 \ 0]^T \tag{8}$$

Then

$$\mathbf{S}_1 = [-2 \ 1]^T \tag{9}$$

Hence the first two terms of the sequence $S^{(1)}$ are {1, -2} and the first two terms of the sequence $S^{(2)}$ are {0, 1}. Subsequent terms in these sequences are to be obtained by the recursion relation

$$S^{(i)}_j = -2 S^{(i)}_{j-1} + S^{(i)}_{j-2}, \qquad j > 1, i = 1, 2 \qquad (10)$$

The application of this recursion leads to the following sequences

| j | $S^{(1)}_j$ | $S^{(2)}_j$ | Ratio |
|---|---|---|---|
| 0 | 1 | 0 | inf |
| 1 | -2 | 1 | -2 |
| 2 | 5 | -2 | -2.5 |
| 3 | -12 | 5 | -2.4 |
| 4 | 29 | -12 | -2.4167 |
| 5 | -70 | 29 | -2.4138 |
| 6 | 169 | -70 | -2.4143 |

The ratio $S^{(1)}_j$ to $S^{(2)}_j$ converges to the largest (absolute value) root, namely $(-1 - \sqrt{2})$, of polynomial p(x) in (6).

**4. Remark**

The eigenvalues of the matrix $\mathbf{R}' = (a\mathbf{I} + b\mathbf{R})$ are $r_i' = (a + br_i)$ with corresponding eigenvectors $[r_i^{m-1}, r_i^{m-2}, \ldots, r_i, 1]^T$, where $\mathbf{I}$ is the *m x m* identity matrix. A different set of sequences may be obtained by replacing matrix $\mathbf{R}$ in (2) by $\mathbf{R}'$, where a, b are suitable integers, chosen so that some other root of p goes to the largest (absolute value) eigenvalue of the $\mathbf{R}'$ without change in its corresponding eigenvector. In the above example, the roots of p are $(-1 - \sqrt{2})$ and $(-1 + \sqrt{2})$. The largest (absolute value) root is $(-1 - \sqrt{2})$. If we take a = 2 and b = 1, we get $\mathbf{R}'$ given by

$$\mathbf{R}' = \begin{bmatrix} 0 & 1 \\ 1 & 2 \end{bmatrix} \qquad (11)$$

The eigenvalues of $\mathbf{R}'$ are $(1 - \sqrt{2})$ and $(1 + \sqrt{2})$. The second eigenvalue now becomes the largest (absolute value) but the eigenvector remains unchanged. Since $(\mathbf{R}' - 2\mathbf{I})$ satisfies the polynomial in (6), therefore, $(\mathbf{R}' - 2\mathbf{I})^2 + 2(\mathbf{R}' - 2\mathbf{I}) - \mathbf{I} = \mathbf{0}$, leading to

$$\mathbf{R}'^2 - 2\mathbf{R}' - \mathbf{I} = \mathbf{0} \qquad (12)$$

We take

$$\mathbf{S}_0 = [1 \ 0]^T \qquad (13)$$

Then
$$\mathbf{S}_1 = [-2 \ 1]^T \qquad (14)$$

Hence the first two terms of the sequence $S^{(1)}$ are $\{1, -2\}$ and the first two terms of the sequence $S^{(2)}$ are $\{0, 1\}$. Subsequent terms in these sequences are to be obtained by the recursion relation

$$S^{(i)}_j = 2 S^{(i)}_{j-1} + S^{(i)}_{j-2}, \qquad j > 1, i = 1, 2 \qquad (15)$$

The application of this recursion leads to the following sequences:

| j | $S^{(1)}_j$ | $S^{(2)}_j$ | Ratio |
|---|---|---|---|
| 0 | 1 | 0 | inf |
| 1 | 0 | 1 | 0 |
| 2 | 1 | 2 | 0.5 |
| 3 | 2 | 5 | 0.4 |
| 4 | 5 | 12 | 0.41667 |
| 5 | 12 | 29 | 0.41379 |
| 6 | 29 | 70 | 0.41423 |
| 7 | 70 | 169 | 0.41420 |

The ratio $S^{(1)}_j$ to $S^{(2)}_j$ converges to the other root, namely $(-1 + \sqrt{2})$, of polynomial p(x) in (6). This suggests that the method can be extended to obtain different sets of integer sequences related to the different roots of a polynomial, except in the case of degenerate roots.

## 5. The case when all roots have the same absolute value

Consider the polynomial $p(x) = x^m - N$. In this case all the *m* roots have the same absolute value. The method as described in sec 2 above will therefore not be applicable. However, one can use the modification described in sec 4 to make the method applicable. The companion matrix, as modified by taking a = 1 and b = 1 is given by

$$\mathbf{R} = \begin{bmatrix} 1 & 0 & 0 & \ldots & 0 & N \\ 1 & 1 & 0 & \ldots & 0 & 0 \\ 0 & 1 & 1 & \ldots & 0 & 0 \\ . & . & . & \ldots & . & . \\ . & . & . & \ldots & 1 & . \\ 0 & 0 & 0 & \ldots & 1 & 1 \end{bmatrix} \qquad (16)$$

As an example take N = 2 and *m* = 3 case. Then $(\mathbf{R} - \mathbf{I})^3 = 2$, so that

$$R^3 - 3R^2 + 3R - I = 0 \tag{17}$$

We take arbitrarily

$$S_0 = [1\ 1\ 0]^T \tag{18}$$

Then

$$S_1 = R S_0 = [1\ 2\ 1]^T \tag{19}$$

and

$$S_2 = R S_1 = [3\ 3\ 3]^T \tag{20}$$

Hence the first three terms of the sequence $S^{(1)}$ are {1, 1, 3}, the first three terms of the sequences $S^{(2)}$ are {1, 2, 3}, and the first three terms of the sequence $S^{(3)}$ are {0, 1, 3}. Subsequent terms in these sequences are obtained by the recursion relation

$$S^{(i)}_j = 3 S^{(i)}_{j-1} - 3 S^{(i)}_{j-2} + 3 S^{(i)}_{j-3}, \quad j > 1,\ i = 1, 2 \tag{21}$$

The application of this recursion leads to the following sequences:

| j | $S^{(1)}_j$ | $S^{(2)}_j$ | $S^{(3)}_j$ | Ratio $S^{(1)}_j / S^{(2)}_j$ | Ratio $S^{(2)}_j / S^{(3)}_j$ |
|---|---|---|---|---|---|
| 0 | 1 | 1 | 0 | 1 | infinity |
| 1 | 1 | 2 | 1 | 0.5 | 2 |
| 2 | 3 | 3 | 3 | 1 | 1 |
| 3 | 9 | 6 | 6 | 1.5 | 1 |
| 4 | 21 | 15 | 12 | 1.400000 | 1.250000 |
| 5 | 45 | 36 | 27 | 1.250000 | 1.333300 |
| 6 | 99 | 81 | 63 | 1.222222 | 1.285714 |
| 7 | 225 | 180 | 144 | 1.250000 | 1.250000 |
| 8 | 513 | 405 | 324 | 1.266667 | 1.250000 |
| 9 | 1161 | 918 | 729 | 1.264706 | 1.259259 |
| 10 | 2619 | 2079 | 1647 | 1.259740 | 1.262295 |
| 11 | 5913 | 4698 | 3726 | 1.258621 | 1.260870 |
| 12 | 13365 | 10611 | 8424 | 1.259542 | 1.259615 |
| 13 | 30213 | 23976 | 19035 | 1.260135 | 1.259574 |
| 14 | 68283 | 54189 | 43011 | 1.260090 | 1.259887 |
| 15 | 154305 | 122472 | 97200 | 1.259921 | 1.260000 |
| 16 | 348705 | 276777 | 219672 | 1.259877 | 1.259956 |
| 17 | 788049 | 625482 | 496449 | 1.259907 | 1.259912 |
| 18 | 1780947 | 1413531 | 1121931 | 1.259928 | 1.259910 |
| 19 | 4024809 | 3194478 | 2535462 | 1.259927 | 1.259919 |
| 20 | 9095733 | 7219287 | 5729940 | 1.259921 | 1.259924 |
| 21 | 20555613 | 16315020 | 12949227 | 1.259920 | 1.259922 |
| 22 | 46454067 | 36870633 | 29264247 | 1.259921 | 1.259921 |
| 23 | 104982561 | 83324700 | 66134880 | 1.259921 | 1.259921 |
| 24 | 237252321 | 188307261 | 149459580 | 1.259921 | 1.259921 |
| 25 | 536171481 | 425559582 | 337766841 | 1.259921 | 1.259921 |

It may be observed that both of the ratios converge to $(2)^{1/3}$ = 1.259921.

## 6. Conclusions

We have obtained a method for finding the real roots of any integer monic polynomial by using recursion relation involving integers only. This method is faster than the conventional methods using floating point arithmetic. The details will be discussed elsewhere.

## 7. Bibiliography